\documentclass[reqno,11pt,a4paper,final]{amsart}
\usepackage[margin=1in]{geometry}
\usepackage{amsmath, enumitem, amssymb, amsthm, amscd, mathtools, bbm, MnSymbol, wasysym}

\usepackage[T1]{fontenc}
\usepackage[utf8]{inputenc}
\usepackage[dvipsnames]{xcolor}

\newcommand{\mres}{\mathbin{\vrule height 1.6ex depth 0pt width
0.13ex\vrule height 0.13ex depth 0pt width 1.3ex}}

\newcommand{\R}{\mathbb{R}}

\newcommand{\norm}[1]{\|#1\|}
\newcommand{\eps}{\epsilon}
	
\newcommand{\loc}{\mathrm{loc}}

\theoremstyle{definition}

\newtheorem{tw}{Theorem}[section]

\newtheorem{defi}[tw]{Definition}
\newtheorem*{dow}{Proof}
\newtheorem{uw}[tw]{Remark}

\newtheorem{lem}[tw]{Lemma}

\newtheorem{prop}[tw]{Proposition}

\DeclareMathOperator{\tr}{tr}
\DeclareMathOperator{\intel}{int}

\DeclareMathOperator{\esssup}{ess\,sup}

\DeclareMathOperator{\diverg}{div}

\DeclareMathOperator{\supp}{supp}

\DeclareMathOperator{\capac}{Cap}

\renewcommand{\epsilon}{\varepsilon}
\renewcommand{\phi}{\varphi}
\renewcommand{\qedhere}{\hfill\blacksquare}

\title{H\MakeLowercase{igher regularity of solutions to elliptic equations on low-dimensional structures}
}
\author{Ł\MakeLowercase{ukasz} C\MakeLowercase{homienia}, M\MakeLowercase{ichał} F\MakeLowercase{abisiak}}
\date{\today}

\address{Institute of Applied Mathematics and Mechanics, University of Warsaw, ul. Banacha 2, 02-097 Warszawa, Poland}
\email{lukasz.chomienia@uw.edu.pl, m.fabisiak@mimuw.edu.pl}

\thanks{\textbf{Acknowledgement.}
We want to thank Professor Anna Zatorska-Goldstein for her support and help, and Professor Piotr Rybka for reading the early version and giving many useful suggestions.\\ 
Both authors are supported by the National Science Center, Poland, by the Grant: 2019/33/B/ST1/00535. 
}



\begin{document}
\maketitle

\begin{abstract}
This paper aims to establish counterparts of fundamental regularity statements for solutions to elliptic equations in the setting of low-dimensional structures such as, for instance, glued manifolds or CW-complexes. 
The main result proves additional Sobolev-type regularity of weak solutions to elliptic problems defined on the low-dimensional structure. Due to the non-standard geometry of a domain, we propose a new approach based on combining suitable extensions of functions supported on the thin structure with uniform bounds for difference quotients. We also derive several important conclusions from that result, namely global continuity of weak solutions and we address the correspondence between the low-dimensional weak problems and the second-order operators widely applied in the framework developed for applications in variational problems.
\end{abstract}



\section{Introduction} 
In recent times, variational problems and partial differential equations in the framework of singular subsets of Euclidean space have attracted the attention of researchers, for example, \cite{first, Bouchitte, man, Rybka}. Interest in such issues is motivated by a variety of possible applications, including theoretical ones like the relaxation of geometrical functionals \cite{omar}, or practical ones like the optimal design issue \cite{lew}, the strength of materials \cite{bolb}, heat transfer on thin conductors \cite{Rybka}, and other topics primarily related to engineering aspects.    

There are several known ways to perform analysis on structures which are singular with respect to the $d$-dimensional Lebesgue measure in the space $\R^d$, such as the fattening method used, for instance, in \cite{fat} or in  \cite[Section 4.2]{Bouchitte} or various measure-theoretic approaches, discussed for example in \cite{first, lew, man}. 
Recently, the latter method has seen many new contributions due to its advantages. Utilizing general measure-theoretic tools allows us to avoid problems related to a low regularity of structures or an occurrence of components of a domain of different dimensions, such as a disc glued to an orthogonal interval. In this paper, we work in a quite modern framework involving singular measures. 

The approach to low-dimensional analysis based on a measure-theoretic framework has well-funded foundations that are positively verified by applications. Among them, we mention the theory of the first-order calculus of variations \cite{first}, second-order variational theory possessing interesting non-local effects \cite{Bouchitte}, as well as the basic theory of elliptic \cite{Rybka} or even parabolic \cite{chom} equations. Currently, the vague direction is to provide correspondence between those theories. The importance of such connections is noticeable, for instance, in \cite{Rybka}, where the authors used direct methods of the calculus of variations to obtain the existence and uniqueness of solutions to the low-dimensional heat transport equation.   

We consider the setting originally developed in \cite{first} and \cite{man}. In \cite{Rybka}, the authors prove the existence of weak solutions to elliptic problems considered on a general class of glued manifolds of potentially different dimensions. This naturally leads to the question of the higher regularity of said solutions. The most fundamental notion of regularity is the higher Sobolev regularity on the components manifolds of the given low-dimensional structure. This kind of regularity is required for further studies but is far from sufficient from the perspective of the low-dimensional structures theory since it does not capture any additional connection between higher-order behaviour on different components. It turns out that a membership of a weak solution in the domain of the second-order operator introduced in \cite{Bouchitte} is the most adequate kind of elliptic regularity in this framework. 
Indeed, the mentioned framework is rich enough to encode the geometry of the structure, provides correspondence between regularity on various components and arises naturally in low-dimensional variational problems.     
   
Firstly, we address the aforementioned problem of upgrading the elliptic regularity in the classical Sobolev sense on components. This is the starting point for further improvements.   
We examine other regularity-related properties of weak solutions. Combining the local regularity result with some facts from the capacity theory, we conclude that weak solutions are continuous. A different critical aspect of our study uses technical properties of the low-dimensional second-order operator to prove the membership of weak solutions in the proper higher-order space of functions.

We are especially interested in the elliptic problems on structures with one or two-dimensional parts embedded in $\R^3$ and throughout the paper we restrict our attention to this setting. Such restriction is motivated by physical applications (see, for example, \cite{lew}). With a small workload, it may be generalised to analogous problems in an arbitrary $d$-dimensional Euclidean space. We also impose other technical restrictions on the class of considered structures to avoid technical complications. Nevertheless, the presented methods likely carry over to more general structures.
 
The following facts concerning an additional regularity of low-dimensional weak solutions to elliptic equations constitute the main results of this paper:
\begin{itemize}
\item On each component manifold $S_i$ of the low-dimensional structure $S$ a weak solution $u$ has the extra regularity $u \in H^2(S_i).$ This is the statement of Theorem \ref{globreg}.$ $\\

\item A weak solution $u$ is globally continuous on the given low-dimensional structure, that is $u \in C(S).$ This fact is precisely expressed in Theorem \ref{ciaglosc_ogolne}.$ $\\

\item If $u$ solves the weak problem, then $u$ is a member of the low-dimensional second-order derivative operator. This is the result of Theorem \ref{nalezenie}.
\end{itemize}
The proof of the result of the first point is quite technical and involves well-known standard regularity results as well as new constructions and methods dedicated to the considered low-dimensional framework. The main obstacle in showing higher regularity lies in defining a proper notion of translation, which is non-trivial due to the geometry of a domain. As we aim to use the difference quotients, it is necessary to address this question. We propose a construction based on a properly chosen sequence of well-behaved extensions of weak solutions. The obtained sequence of approximations shares two crucial properties. Firstly, we can extend each approximation from the considered structure to the whole Euclidean space and obtain a sufficiently regular function. Secondly, the global behaviour of extensions is controlled in terms of the original function posed on the low-dimensional structure. 
As a result, we obtain that the functional space is closed under taking such generalised difference quotients. This allows us to prove a higher regularity of solutions by controlling norms of difference quotients. 
Latter results are implications or refinements of the main theorem. We establish the global continuity of low-dimensional weak solutions by utilising some facts from the Sobolev capacity theory and conclusions provided by our main result. An interesting fact is that the continuity on the whole structure depends on the dimensions of the components.
Furthermore, still relying on the main theorem, we provide the membership of a weak solution to the domain of the measure-related second-order operator.
In all of those results, the key obstacle is related to the fact that measures associated with the examined structure are singular with respect to the standard Lebesgue measure in $\R^3$.

Let us briefly discuss some other potentially related results and highlight the differences. Equipped with the naturally induced metric, the low-dimensional structure $S$ is a metric measure space. In the category of metric measure spaces, there are several known candidates for generalisations of the Euclidean gradient, amongst them two are probably best known: the upper gradient discussed in \cite{ambro,hei} and the Cheeger gradient, see \cite{cheeg}. It turns out that interpreting the low-dimensional setting as metric measure spaces, the upper gradient corresponds to $|\nabla_{S_i}|,$ where $\nabla_{S_i}$ is a component of the gradient tangent to $S_i$. In other words, such an approach loses a significant amount of information required by our goals. Therefore, we start with the $\mu$-related gradient, which is the most natural generalisation of a gradient tangent to a manifold. It should be mentioned that our regularity results are strongly based on the form of the gradient and the fact that it represents the classical gradient locally. Potential generalisations of the obtained regularity results to a less geometrical setting seem difficult.   

The paper is organized as follows. In Section 2, we formulate basic definitions and facts, amongst which we introduce the notion of a low-dimensional structure, adequate spaces of functions and discuss their properties. Section 3 is devoted to the main theorem and some closely related facts. Both Section 4 and Section 5 address important implications of the main result; in Section 4 we prove the global continuity of solutions under the additional assumption that all components are of the same dimension, and in Section 5 we build the correspondence between weak solutions and the second-order spaces. Since their definition is highly non-trivial and, at the same time, not crucial for the rest of this paper, we recall the second-order framework in the Appendix.   

\section{Preliminaries}

There are many ways of introducing a tangent bundle of a measure and the Radon measure-based first-order space of functions. The difference between various notions of the tangent structure to a measure has been studied in \cite{frag}. As those notions are not equivalent in general, it is necessary to point out throughout this paper we stick to the first-order framework proposed in \cite{first} and \cite{man}. We opt to introduce the counterpart of the Sobolev spaces in this way because the procedure can be naturally mimicked to define the second-order measure-related functional spaces (see the Appendix). We also recall the essential definitions regarding the elliptic equations, which are considered in a form proposed in \cite{Rybka}.
Let us point out that the main reason for restricting our attention to the case of structures embedded in $\R^3$ is the fact that part of the applied here results related to the second-order framework, like Theorem \ref{domk} is known only in such setting.  

Assume that $\mu$ is a positive Radon measure in $\R^3.$ The $\mu$-related counterpart of the traditional Lebesgue space is the most fundamental function space we use.
For an open $\Omega \subset \R^n$ and $p \in [1,\infty],$ we define the Lebesgue space $L^p_{\mu}$ as a subspace of the space of $\mu$-measurable functions such that the norm 
$$
\displaystyle \norm{\cdot}_{L^p_{\mu}}:= \begin{cases}
\displaystyle \left( \int_{\Omega} |f|^p d\mu\right)^{\frac{1}{p}} &\text{ for } 1\leqslant p < \infty, \\
\displaystyle \esssup_{\mu}|f|    &\text{ for } p = \infty
\end{cases}
$$ 
is finite. In further perspective, the measure $\mu$ will be understood to correspond to a kind of geometrical object and the notion of tangent structure will be especially important for us. Nevertheless, we can define a concept of the tangent structure in a general instance of an arbitrary Radon measure. Let 
$$
\mathcal{N}_{\mu}:=\{v \in C^{\infty}_c(\R^3): v=0 \text{ in } \supp \mu \}
$$ 
and 
$$
\mathfrak{N}_{\mu}:= \{ w \in C^{\infty}_c(\R^3;\R^3): w=\nabla v \text{ in } \supp \mu \text{ for some } v \in \mathcal{N}_{\mu} \}.
$$
We introduce the multifunction $N_{\mu}: \R^3 \to \mathcal{P}(\R^3)$ (where $\mathcal{P}(X)$ denotes the power set of $X$):
$$
N_{\mu} (x):= \{w(x)\in \R^3: w \in \mathfrak{N}_{\mu}  \}.
$$
At $\mu$-a.e. point $x\in \R^3$, the space $T_{\mu}(x)$ tangent to the measure at $x \in \R^3$ is defined as 
$$
T_{\mu}(x):= \{ u(x) \in \R^3: u(x)\perp N_{\mu} (x)  \}.
$$
The symbol $\perp$ stands for orthogonality in the sense of the standard scalar product in $\R^3.$  Furthermore, by $T_{\mu}^{\perp}(x)$ we denote the $\R^3$-orthogonal completion of $T_{\mu}(x).$ 

With the tangent space at hand, we can project onto it. In particular, we can determine the $\mu$-essential component of the gradient. To this end, denote the orthogonal projection onto the tangent space as $P_\mu: \R^3 \to T_\mu$. The tangent gradient of a function $u \in C^\infty_c(\R^3)$ is defined as
$$
\nabla_{\mu} u(x):= P_{\mu} \nabla u.
$$
The Sobolev space $H^1_{\mu}$ is defined as a completion of the space $C^{\infty}_c(\R^3)$ with respect to the Sobolev-type norm 
$$
\norm{\cdot}_{\mu}:= \left(\norm{\cdot}^2_{L^2_{\mu}} + \norm{\nabla_{\mu} \cdot}^2_{L^2_{\mu}}\right)^{\frac{1}{2}}.
$$	
We define the zero-mean subspaces: 
$$
\mathring{L^2_{\mu}}:=L^2_{\mu} \cap \left\{u\in L^2_{\mu}: \int_{\Omega}ud\mu=0 \right\}\; \text{  and  }\; \mathring{H^1_{\mu}}:=H^1_{\mu} \cap \mathring{L^2_{\mu}}.
$$

We restrict our interest to the class of measures possessing a geometrical background. 
\begin{defi}\label{class}(Class $\mathcal{S}$ of low-dimensional structures)\\
Let $\Omega \subset \R^3$ be a non-empty, open and connected set. For a fixed $m \in \mathbb{N}$ and $0\leqslant k \leqslant m+1,$ let $S_i,\; 1 \leqslant i \leqslant m$ be a:
\begin{itemize}
    \item $2$--dimensional closed manifold with boundary for $i\leqslant k$;
    \item $1$--dimensional closed manifold with boundary for $i > k$.
\end{itemize}
Assume further that for each $1 \leqslant i,j \leqslant m, \; i \neq j$:
\begin{itemize}
\item[a)] $\partial \Omega \cap S_i = \partial S_i$;
\item[b)] $S_i$ is transversal to $S_j$;
\item[c)] for any $1\leqslant i,j,k \leqslant m, \; S_i\cap S_j \cap S_k = \emptyset.$
\end{itemize}   
With each $S_i$ we associate the pair $(\mathcal{H}^{\dim S_i}\mres_{S_i}, \theta_i),$ where $\mathcal{H}^{\dim S_i}\mres_{S_i}$ is the $\dim S_i$-dimensional Hausdorff measure restricted to $S_i$ and $\theta_i \in L^{\infty}(S_i)$, where $\theta_i \geqslant c > 0,$ for some constant $c.$
We say that a positive Radon measure $\mu$ belongs to the class $\mathcal{S}$ if it is of the form $$\mu = \sum_{i=1}^m \theta_i \mathcal{H}^{\dim S_i}\mres_{S_i}.$$
\qed
\end{defi}

\begin{uw}
Let us briefly explain the restrictions imposed on the $\mathcal{S}$ class. First of all, we should understand that each structure $S$ is dependent on the choice of the region $\Omega,$ maybe not in a topological, but at least in a geometrical sense. This is forced by the assumption a). The transversality condition b) excludes those kinds of domains that do not really need the presented framework and similar results can be proved with the help of standard methods. Finally, condition c) disallows for more than two component manifolds to intersect in a single point. This assumption is introduced to simplify technicalities, and we believe it can be relaxed for the price of longer computations. Throughout the paper, we call measures belonging to the class $\mathcal{S}$ the low-dimensional structures. We also use the same name to refer to sets $S=\bigcup_{i=1}^mS_i$ on which the singular measures are supported. This ambiguity does not produce confusion; the type of object we refer to should be clear from a context.
\end{uw}

In order to introduce an appropriate elliptic setting, we begin with a notion of a suitable relaxation of a given matrix of coefficients.
\begin{defi}\label{relax}(Relaxed matrix of coefficients; \cite[Proposition 3.1]{Rybka} )$ $\\
Let $B \in \left(L^{\infty}_{\mu}\right)^{3 \times 3}$ be such that for $\mu$-almost every $x \in \R^3$:
\begin{itemize}
\item[a)] $B(x)$ is a symmetric matrix;
\item[b)] $T_{\mu}(x) \subset  \text{Im}B(x)$;
\item[c)] $(B(x)\xi,\xi)\geqslant \lambda |\xi|^2$ for all $\xi \in \text{Im}B(x)$ and some $\lambda>0$. 
\end{itemize}
We define the relaxation $B_{\mu}$ of the matrix $B$ by the formula 
$$
B_{\mu}(x) := B(x) - \sum_{i=1}^l \frac{B(x)e_i(x)\otimes B(x)e_i(x)}{(B(x)e_i(x),e_i(x))},
$$ 
where $l = 1,2,$ $e_i(x),$ $1 \leqslant i \leqslant l$ are linearly independent, $\mu$-measurable, span $T_{\mu}(x)^{\perp} \cap ImB(x)$ and for $1 \leqslant i,j \leqslant l $, $(B(x)e_i(x),e_j(x)) = \delta_{ij},$ where $\delta_{ij}$ is the Kronecker's delta. 
\qed
\begin{uw}
The matrix $B_{\mu}$ is indeed a relaxation of $B$ in the proper sense. This fact is proved in \cite[Proposition 3.1]{Rybka}.
\end{uw}
\end{defi}

We now refer to the low-dimensional weak equations and the theory first presented in \cite{Rybka}. 

\begin{defi}[Weak formulation of Neumann problem, \cite{Rybka}]\label{defin
}$ $\\
Let $f \in L^2_{\mu}$ and $\int_{\Omega} fg d\mu = 0$ for all $g \in \ker \nabla_{\mu}.$ We say that $u \in H^1_{\mu},$ such that $\int_{\Omega} ug d\mu = 0$ for all $g \in \ker \nabla_{\mu},$ is a weak solution to the elliptic Neumann problem 

$$\diverg(B\nabla u)=f,\quad B\nabla_{\mu} u \cdot \nu\lvert_{\partial S}=0$$ 
provided that

\begin{gather}\label{weak}
\int_{\Omega}B_{\mu}\nabla_{\mu} u \cdot \nabla_{\mu}\phi d\mu = \int_{\Omega} f\phi d\mu,
\end{gather}
holds for all $\phi \in C^1_c(\R^3).$  \qed
\begin{uw}\label{H1 jako testowe}
By the definition of the space $H^1_{\mu}$, smooth functions $C^{\infty}_c(\R^3)$ are dense in $H^1_{\mu}.$ This means that the class of test functions in the definition above may be extended to $H^1_{\mu}.$\\ 
Moreover, by the fundamental theorem of the calculus of variations it can be checked that condition $\int_{\Omega}u \phi d\mu =0 $ for all $\phi\in \ker\nabla_{\mu}$ implies $\int_{\Omega} u d\mu=0.$\\
As it is stated in paper \cite{Rybka} it is possible to give a meaning to the zero Neumann condition using the theory of measures with divergences, see \cite{Chen}.
\end{uw} 
\end{defi}

Combining a useful characterisation of the space $H^1_{\mu}$ established in \cite{Bouchitte} with the fact that multiplying a measure by a bounded and separated from zero density does not change the tangent structure, we obtain the following lemma.
\begin{lem}[Characterisation of $H^1_{\mu};$ {\cite[Lemma 2.2]{Bouchitte}}]\label{kombo char}$ $\\
Let $\mu \in \mathcal{S}.$ The tangent structure on the component manifold $S_i$ is denoted by $T_{S_i}$ and the classical tangent gradient on the component $S_i$ is denoted by $\nabla_{S_i}.$ Then:
\begin{itemize}
\item[a)] $\displaystyle T_{\mu}(x) = \sum_{i=1}^mT_{S_i}(x)$ for $\mu$-a.e. $x$;
\item[b)] if $\displaystyle u \in H^1_{\mu},$ then for $1 \leqslant i \leqslant m$ we have $\displaystyle u\lvert_{S_i} \in H^1(S_i)$;
\item[c)] if $\displaystyle u \in H^1_{\mu},$ then for $1 \leqslant i \leqslant m$ we have $(\displaystyle\nabla_{\mu}u)\lvert_{S_i} = \nabla_{S_i}u$;
\item[d)] if $\displaystyle\phi \in C^{\infty}(\R^3),\; \phi=0 \text{ on } \bigcup_{j\neq i}^mS_j  \text{ and } u \in L^2_{\mu}\cap H^1(S_i),$ then $\displaystyle\phi u \in H^1_{\mu}.$  
\end{itemize}
\end{lem}

On an arbitrary low-dimensional structure $S \in \mathcal{S},$ the classical Poincar{\'e} inequality in the space $H^1_{\mu}$ is invalid. This is closely related to the fact that the intersection set of components of different dimensions is of zero $W^{1,2}$-Sobolev capacity on the higher dimensional component. This fact is discussed in \cite[Section 2.2]{Rybka}. In the same paper, the authors establish the weaker form of the Poincar{\'e} inequality. We evoke it now for further applications. 

\begin{defi}\label{slabypoincare}(Weak Poincar{\'e} inequality; \cite[Theorem 2.1]{Rybka})\\
 For the low-dimensional structure with $m$ component manifolds we group the set $\{1,...,m\}$ into subfamilies $I_1,...,I_d,$ that is $\left\{1,...,m\right\}=I_1\cup...\cup I_d,$ where sets $I_k,\; 1 \leqslant k \leqslant d$ are pairwise disjoint. Besides that, we impose the following assumptions: 
        \begin{itemize}
		\item[a)] For any $1 \leqslant k \leqslant d$, the characteristic function $\chi_{\left\{\bigcup_{i\in I_k}S_i\right\}}$ of a sum of component manifolds with indexes in $I_k$ is an                       element of the kernel of the $\nabla_{\mu}$ operator:  
            		\begin{equation*}
            		      \chi_{\left\{\bigcup_{i\in I_k}S_i\right\}} \in \ker \nabla_{\mu}.
            		\end{equation*}
		\item[b)] Each $I_k,\;1 \leqslant k \leqslant d$ is maximal in a sense that: if $\alpha \in I_k$ and $\widetilde{I}_k:= I_k \setminus \{\alpha\},$ then 
            		\begin{equation*}
            		      \chi_{\left\{\bigcup_{i\in \widetilde{I}_k}S_i\right\}} \notin \ker \nabla_{\mu}.
            		\end{equation*}
	\end{itemize}  
For each $I_k$ we define the projection 
    \begin{equation*}
    P_ku:=\chi_{\left\{\bigcup_{i\in I_k}S_i\right\}} \strokedint_{\bigcup_{i\in I_k}S_i} u d \mu,
    \end{equation*}
where $\strokedint_Xfd\mu := \frac{1}{\mu(X)}\int_Xfd\mu.$
Then for any $\mu \in \mathcal{S}$ there exist positive constants $C_i,\; 1 \leqslant i \leqslant d$ such that for each $1 \leqslant k \leqslant d$ and all $u \in H^1_{\mu}$ the weak Poincar{\'e} inequality
	\begin{equation}\label{weakpoincare}
	   \sum_{j\in I_k} \int_{\Omega}|u-P_ku|^2d \mu_j \leqslant C_k \sum_{j\in I_k}\int_{\Omega}|\nabla_{\mu}u|^2d\mu_j. 
	\end{equation}
	  is satisfied. Here we denote $\mu_j:=\mu\mres_{S_j}.$
  \qed
\end{defi}

Throughout the paper, $C$ is a generic, positive constant and may change from line to line.
\section{Sobolev-type regularity of weak solutions}

One of the standard ways of showing extra Sobolev regularity in the Euclidean setting is the method of difference quotients - if they are uniformly bounded, they justify the existence of a weak gradient. In a classical setting (see for example \cite[Section 6.3]{Evans}), one obtains uniform bound on difference quotients of weak solutions, thus improving their regularity. We aim to proceed analogously on low-dimensional structures. However, the difficulties in such an approach are readily visible - in general, we cannot easily shift the function supported on the structures we work with in the directions orthogonal to its components. We circumvent this issue in the following way:
\begin{itemize}
    \item For a given $u \in H^1_{\mu}$, we take a sequence of smooth functions $ C^\infty_c(\R^3) \ni \phi_n \xrightarrow{H^1_\mu} u$, which we can arbitrarily shift in any direction; 
    \item As we shall see, we need further control of $\phi_n$, and therefore we introduce a special modification of $\phi_n$, which we denote by $\widehat{\phi_n}$. Its key property is that its global behaviour is determined by the restriction to a low-dimensional structure;
    \item Our modification still obeys the $H^1_\mu$ convergence and its shifts are well-defined. The limit $\lim\limits_{n\rightarrow \infty} \widehat{\phi_n}(\cdot + h)$ plays the role of a generalized translation of $u$. 
\end{itemize}

\begin{uw}
The third point above is not as immediate as it may seem, since for an arbitrary approximating sequence $\phi_n \in C^\infty_c(\R^3)$, justifying the membership $u \in H^1_\mu$, it may happen that $\phi_n (\cdot + h)$ no longer converges in $H^1_\mu$. This is one of the fundamental reasons why we introduce a specific form of the extension. 
\end{uw}


The proof of our main theorem -- $H^2$-type regularity on each component manifold -- consists of two parts: showing that the functional space $H^1_{\mu}$ is closed with respect to the aforementioned generalised translation and secondly, establishing uniform bounds on generalized difference quotients. 

In what follows, we restrict our attention to the structures consisting only of "straightened out" components. While such structures serve as generic models of various types of intersections, we later show that the general result follows from the one we obtain in the model cases. We address this issue at the end of this section.
We consider the following structure:\\ 
Let $S = S_1 \cup S_2$, where 
\begin{equation}\label{dyski}
\begin{aligned}
S_1 = \{ (x,y,0) \in \R^3 \ : \ x^2 + y^2 \leqslant 1 \}, \quad S_2 =  \{ (0,y,z) \in \R^3 \ : \ y^2 + z^2 \leqslant 1 \}.
\end{aligned}
\end{equation}
We introduce two Hausdorff measures associated with $S$: 
\begin{equation}\label{miary}
\begin{aligned}
\mu := \mathcal{H}^2\mres_{S_1} +\mathcal{H}^2\mres_{S_2}, \quad \widetilde{\mu} := \theta_1\mathcal{H}^2\mres_{S_1}+\theta_2\mathcal{H}^2\mres_{S_2},\quad \theta_i\in C^{\infty}(S_i)\text{ for } i=1,2.
\end{aligned}
\end{equation}
Denote the intersection set by 
$$
\Sigma:=S^1\cap S^2 = \{(0,y,0) \ : \ |y|\leqslant 1 \}
$$ and the restriction of $u$ to a single disc as $u_i:=  u\lvert_{S_i}$. We are ready to state the first result.

\begin{prop}\label{slady}
Let $u \in H^1_{\mu}.$ Then $\tr^{\Sigma}u_1=\tr^{\Sigma}u_2$, where $\tr^{\Sigma}$ denotes the usual trace on $\Sigma$.
\end{prop}
\begin{dow}
Let $\alpha_n \in C^{\infty}(\R^3),\; \alpha_n \rightarrow u$ in $H^1_\mu$ and let $\alpha_n^i:= \alpha_n \lvert_{S_i}, \; i=1,2$. Since $\alpha_n$ is smooth, we know that $\tr^{\Sigma}\alpha^i_n = \alpha^i_n \lvert_{\Sigma}$. Since $\alpha_n^1\lvert_{\Sigma} = \alpha_n^2\lvert_{\Sigma}$, it follows that $\tr^{\Sigma}\alpha^1_n = \tr^{\Sigma}\alpha^2_n$. Since $\alpha^i_n$ converges to $u_i$ in $H^1(S_i)$, the continuity of the trace operator implies 
$$
\tr^{\Sigma}\alpha^i_n \xrightarrow{L^2(\Sigma)}\tr^\Sigma u_i, \; i=1,2.
$$
However, $\tr^{\Sigma}\alpha^1_n = \tr^\Sigma \alpha^2_n$ and thus $\tr^{\Sigma}u_1 = \tr^{\Sigma}u_2$.
$\qedhere$
\end{dow}

Let us now consider the low-dimensional elliptic problem for $S$ previously studied in \cite{Rybka}. For a given $f \in \mathring{L^2_{\mu}},$ a solution $u \in \mathring{H^1_{\mu}}$ satisfies the equality 
\begin{equation}\label{cieplo2}
\int_{\Omega}B_{\widetilde{\mu}} \nabla_{\widetilde{\mu}}u \cdot \nabla_{\widetilde{\mu}} \phi d\widetilde{\mu} = \int_{\Omega}f \phi d\widetilde{\mu}
\end{equation}
for any $\phi \in C^{\infty}(\R^3).$ Existence of the unique solution to this problem is the main result of paper \cite{Rybka}. By Proposition \ref{kombo char} point c), we may include the densities in the operator $B_{\mu}$ and change $\widetilde{\mu}$-related gradients to $\nabla_{\mu},$ 
\begin{equation}\label{cieplo3}
\int_{\Omega}\widetilde{B_{\mu}} \nabla_{\mu}u \cdot \nabla_{\mu} \phi\, d\mu = \int_{\Omega}\widetilde{f} \phi\, d\mu.
\end{equation}
We now apply the strategy outlined before to the solution $u$.\\

\noindent\textbf{Step 1:} $\mathbf{\tr^{\Sigma} u \in H^1_{\loc}(\Sigma)}$ \\
We begin by showing extra differentiability in the direction of the variable $y$. Let us define 
$$
u^{h,y}:= u(\cdot + he_y), \text{ where } e_y=(0,1,0).
$$
Such translation is a well-defined function, at least on any low-dimensional structure $S'$ satisfying $\intel S' \subset \subset \intel S$ (in the inherited topology on $S;$ we refer to $S'$ as a substructure of $S$) and for sufficiently small $h>0.$ We abuse the notation a little and by $\mu$ denote the restriction $\mu\mres_{S'}$; it is justified since all of the further reasoning is carried out locally. Taking any $\phi_n\in C^{\infty}(\R^3)$ satisfying $\phi_n \rightarrow u$ in $H^1_\mu$, we see that 
$$
\phi_n^{h,y} \xrightarrow[n\to \infty]{H^1_{\mu}} u^{h,y},
$$ 
and thus $u^{h,y} \in H^1_{\mu}.$ 

We put 
$$
D^h_yu:= \frac{1}{h}(u^{h,y}-u) \in H^1_\mu.
$$ Let $\xi \in C^{\infty}_c(\R^3)$ be a smooth function which will be specified later. We have $\xi^2D^h_yu \in H^1_{\mu}$ and further 
$\phi^y:=-D^{-h}_y(\xi^2D^h_yu) \in H^1_{\mu}.$ This means that $\phi^y$ can be used as a test function in \eqref{cieplo3}.

Since $\Sigma$ is $\mu$-negligible, we are able to rewrite 
$$
\int_{\Omega}\widetilde{B_{\mu}}\nabla_{\mu} u \cdot \nabla_{\mu} \phi^y d\mu = \int_{S_1}\widetilde{B_{\mu}}\nabla_{S_1} u \cdot \nabla_{S_1} \phi^y dS_1 + \int_{S_2}\widetilde{B_{\mu}}\nabla_{S_2} u \cdot \nabla_{S_2} \phi^y dS_2,$$ 
and the right-hand side of the equation \eqref{cieplo3} can be presented as 
$$
\int_{\Omega} \widetilde{f} \phi^y d \mu = \int_{S_1}\widetilde{f} \phi^y dS_1+ \int_{S_2}\widetilde{f} \phi^y dS_2.$$ 

Having the equality in the expanded form as above, we might apply the standard method of showing higher regularity of solutions by establishing a uniform bound on difference quotients (for example, see Section 6.3.1 in \cite{Evans}). 

Recall that we assumed $\mu$ to be a restriction of the measure $\mathcal{H}^1\mres_{S_1}+ \mathcal{H}^1\mres_{S_2}$ to a fixed subset $S' \subset \subset S.$ This provides that for small enough $h>0$ we have $\norm{D^h_y \nabla_{\mu}u}_{L^2_{\mu}}\leqslant C,$ which further implies that $\norm{\partial_y \nabla_{\mu} u}_{L^2_{\mu}} \leqslant C,$ where the constant $C>0$ is independent of $h.$ Using the explicit form of $\nabla_\mu$, we have 
$$\norm{\partial_y^2u}_{L^2_{\mu}}\leqslant C \text{ and } \norm{\partial_y \partial_x u}_{L^2_{\mu}}\leqslant C.$$ 

Now, we make use of the fact that in a region separated from the junction set $\Sigma$, the equation is reduced to the standard case, and the local regularity is known. Define 
$$
S_1^+:= \{(x,y,0)\in S_1: x>0\}, \quad S_1^-:= \{ (x,y,0) \in S_1: x<0\}.
$$ 
Let $V \subset\subset S_1^+;$ in particular, notice that $\inf\{x: \exists y, (x,y,0)\in V\}\geqslant \alpha>0.$ Taking for $\xi$ properly chosen (as in the classical proof) function supported in $V,$ we conclude by the standard elliptic regularity theory that $u \in H^2_{\loc}(V).$ 

Moreover, the $H^2$-regularity implies almost everywhere in $V$ the symmetry of the second derivatives: $\partial_y\partial_xu=\partial_x\partial_yu.$ By the fact that the set $V$ was chosen arbitrary, we deduce that $\partial_y\partial_xu=\partial_x\partial_yu$ is valid a.e. in $S_1^+.$ Proceeding analogously on $S_1^-$ we get $\partial_y\partial_xu=\partial_x\partial_yu$ a.e. on $S_1^-,$ so the equality $\partial_y\partial_xu=\partial_x\partial_yu$ is true a.e. on $S_1.$ By the fact that $\partial_y\partial_xu \in L^2_{\loc}(S_1)$ and $\partial_y\partial_xu=\partial_x\partial_yu$ a.e. on $S_1$ we deduce that $\partial_x\partial_y u \in L^2_{\loc}(S_1).$ We will check that $\partial_x\partial_y u$ is actually a weak $\partial_x$-derivative of $\partial_yu$ on $S_1.$ 

Let 
$$
\mathcal{A}:= \{\phi \in C^{\infty}_c(\R^3): \phi\lvert_{S_1}\in C^{\infty}_c(S_1)\}. 
$$
For any $\phi \in \mathcal{A}$ it follows that 
\begin{align*}
\int_{S_1}\partial_x\partial_yu \phi dS_1 &= \int_{S_1}\partial_y\partial_xu \phi dS_1 = - \int_{S_1}\partial_xu \partial_y\phi dS_1 \\
&= \int_{S_1}u \partial_x\partial_y\phi dS_1 = \int_{S_1}u \partial_y\partial_x\phi dS_1 = - \int_{S_1}\partial_yu\partial_x\phi dS_1,
\end{align*}
thus $\partial_x\partial_yu = \partial_x(\partial_yu),$ so indeed $\partial_x\partial_yu$ is a weak derivative of $\partial_yu.$ 

The result implies that $\partial_yu \in H^1(S_1)$, and so $\Sigma$-trace of $\partial_yu$ is well-defined. It follows that $\tr^{\Sigma}\partial_yu \in L^2(\Sigma).$

Notice that the standard mollification argument shows that smooth functions $C^{\infty}(\R^2)$ are dense in the subspace 
$$
\{u\in H^1(S_1): \partial_yu \in H^1(S_1)\}
$$ inherited with the norm 
$$
\left(\norm{u}_{H^1(S_1)}^2+\norm{\partial_y u}^2_{H^1(S_1)}\right)^{\frac{1}{2}}.
$$ This fact provides the existence of a smooth sequence $\phi_n$ approximating $u$ in the above norm, and further, we can derive the commutation $\tr^{\Sigma}\partial_yu=\partial_y\tr^{\Sigma}u.$ As a conclusion this gives $\partial_y\tr^{\Sigma}u \in L^2(\Sigma).$ A verification that $\partial_y\tr^{\Sigma}u$ is a weak derivative of $\tr^{\Sigma}u$ is a byproduct of the proof of the above commutation. Finally, we obtain 
$$
\tr^{\Sigma}u \in H^1(\Sigma).
$$

\noindent\textbf{Step 2: Construction of the extension}\\
Now we introduce the special form of an extension of a function supported on a low-dimensional structure to the whole $\R^3$. Before proceeding with a rather technical construction, we propose an informal description: for a given point $(x,y,z)\in \R^3$ (not necessarily $(x,y,z)\in S$), the extension can be expressed as 
$$
\widetilde{u}(x,y,z) := u(0,y,z) - (\tr^\Sigma u)(0,y,0) + u(x,y,0).
$$
 
By the previous step, we know that $\tr^{\Sigma}u \in H^1(\Sigma).$ Let us recall that the classical trace theory provides the trace operator 
$$ \tr^{\Sigma}:H^1(S_1) \to H^{1/2}(\Sigma)$$ 
which is surjective and $\tr^{\Sigma}(H^{3/2}(S_1))= H^1(\Sigma).$ Thus there exists $v \in H^{3/2}(S_1)$ such that $\tr^{\Sigma}v = \tr^{\Sigma}u.$ It follows that $\tr^{\Sigma}(v-u) = \tr^{\Sigma}v - \tr^{\Sigma}u  = 0.$ Moreover, we have $v-u \in H^1(S^1)$.

Observe that proceeding as in  Step 1, we obtain that for any set $W\subset \subset S_1$ separated from the intersection $\Sigma$ it holds that $u \in H^2(W)$. Let us choose small, fixed $h>0$ and assume that $\Sigma-(h,0,0) \subset W$ for some fixed $W$ as above.

Now, let $\alpha_n, \beta_n \in C^{\infty}(\R^2)$ be sequences converging in the $H^1(S_1)$-norm to $v$ and $v-u$, respectively. Besides that, we demand from the sequence $\alpha_n$ to converge to $v$ in $H^{3/2}(S_1)$ and from $\beta_n$ to converge to $v-u$ in $H^{3/2}(W)$; the latter can be constructed using the diagonal argument. Without the loss of generality, we can assume $\beta_n\lvert_{\Sigma}=0.$ The extension of such sequences to the $\R^2$ space is a standard fact; first, we judge the existence of sequences converging on $S_1,$ then each sequence term can be smoothly extended to the whole $\R^2$ space. By the continuity of the trace operator, we have
$$
\tr^{\Sigma}\alpha_n \xrightarrow{H^1(\Sigma)}\tr^{\Sigma}v, \quad \tr^{\Sigma}\beta_n \xrightarrow{H^1(\Sigma)}\tr^{\Sigma}(v-u) = 0.
$$
We express the function $u$ as $$u=v-(v-u)$$ and we put $\gamma_n:= \alpha_n - \beta_n \in C^{\infty}(\R^2).$
Immediately, we see that $\tr^{\Sigma}\gamma_n = \tr^{\Sigma}(\alpha_n - \beta_n) = \tr^{\Sigma}\alpha_n,$ and further 
$$
\gamma_n \xrightarrow{H^1(S_1)} v- (v-u) = u \text{ and } \tr^{\Sigma}\gamma_n \xrightarrow{H^1(\Sigma)}\tr^{\Sigma}v = \tr^{\Sigma}u.
$$

Recall that $\alpha_n$ is defined only on $S^1$, which lies in $xy$ plane. We introduce the $\R^3$-extension of $\alpha_n$ by the formula $\widetilde{\alpha_n}(x,y,z):= \alpha_n(x,y).$ We collect some readily seen properties of this extension: $\widetilde{\alpha_n} \in C^{\infty}(\R^3)\cap H^1_\mu,$ 
$$
\widetilde{\alpha_n}\lvert_{S} = \begin{cases}
    \alpha_n &\text{ on } S_1,\\
    \tr^{\Sigma}\alpha_n &\text{ on } S_2
\end{cases}
\qquad \text{and}
\qquad 
\widetilde{\alpha_n} \xrightarrow{H^1_\mu} \begin{cases}
    v &\text{ on } S_1,\\
    \tr^{\Sigma}v &\text{ on } S_2
\end{cases}
= 
\begin{cases}
    v &\text{ on } S_1,\\
    \tr^{\Sigma}u &\text{ on } S_2,
\end{cases}
$$
with the last expression belonging to $H^1_\mu$ (by the completeness of this space).

Analogously we extend $\beta_n$ to $\widetilde{\beta_n}$ and we introduce $\widetilde{\gamma_n}:= \widetilde{\alpha_n}-\widetilde{\beta_n},$ which is exactly the extension of $\gamma_n$ one would obtain by following the same process as for $\alpha_n$ and $\beta_n$. Similarly as before, the function $\widetilde{\gamma_n}$ has the following properties: ${\widetilde{\gamma_n} \in C^{\infty}(\R^3)\cap H^1_{\mu},}$	 
$$
\widetilde{\gamma_n}\lvert_{S}=\begin{cases}
\gamma_n &\text{ on } S_1,\\
\tr^{\Sigma}\alpha_n &\text{ on } S_2
\end{cases}
\quad \text{and}
\quad
\widetilde{\gamma_n} \xrightarrow{H^1_{\mu}} \begin{cases}
u &\text{ on } S_1,\\
\tr^{\Sigma}u &\text{ on } S_2
\end{cases}\in H^1_{\mu}.
$$  	
Repeating analogous construction on the component $S_2$, we obtain
$\widetilde{\delta_n}\in C^{\infty}(\R^3)\cap H^1_\mu$ such that
$$
\widetilde{\delta_n}\lvert_{S}=\begin{cases}
\tr^{\Sigma}\alpha_n' & \text{ on } S_1,\\
\delta_n &\text{ on } S_2
\end{cases}
\quad \text{and}
\quad
\widetilde{\delta_n} \xrightarrow{H^1_{\mu}} \begin{cases}
\tr^{\Sigma}u &\text{ on } S_1,\\
u, &\text{ on } S_2
\end{cases} \in H^1_{\mu},
$$
where $\alpha'_n$ is defined just as $\alpha_n$ was defined on $S_1$.

Having extended the function $u$, we proceed with extending the trace in a similar fashion. As $\tr^{\Sigma}u \in H^1(\Sigma)$ we might find a sequence $\rho_n \in C^{\infty}(\R)$ converging to $\tr^\Sigma u$ in $H^1(\Sigma)$. Define $\widetilde{\rho_n}(x,y,z):= \rho_n(y),\; \widetilde{\rho_n} \in C^{\infty}(\R^3)\cap H^1_\mu.$ The function $\widetilde{\rho_n}$ shares the following properties: 
$$
\widetilde{\rho_n}\lvert_{S}=\begin{cases}
\rho_n & \text{ on } S_1,\\
\rho_n & \text{ on } S_2
\end{cases}
\quad \text{and}
\quad
\widetilde{\rho_n} \xrightarrow{H^1_{\mu}}\begin{cases}
\tr^{\Sigma}u & \text{ on } S_1,\\
\tr^{\Sigma}u & \text{ on } S_2
\end{cases}\in H^1_{\mu}.$$
With all the needed extensions in hand, we define the sequence $\Gamma_n \in C^{\infty}(\R^3):$ 
\begin{equation}\label{extension}
\Gamma_n(x,y,z):= \widetilde{\gamma_n}(x,y,z)+\widetilde{\delta_n}(x,y,z)-\widetilde{\rho_n}(x,y,z).
\end{equation}

\noindent We have $\Gamma_n \in H^1_{\mu}$ and 
$$
\Gamma_n\lvert_{S}= \begin{cases}
\gamma_n+ \tr^{\Sigma}\alpha_n'-\rho_n &\text{ on } S_1,\\
\tr^{\Sigma}\alpha_n+\delta_n-\rho_n, &\text{ on } S_2.
\end{cases}
$$

\noindent Moreover, collecting all the limits of $\widetilde{\gamma_n}, \widetilde{\delta_n}$ and $\widetilde{\rho_n}$ we arrive at
$$
\Gamma_n \xrightarrow{H^1_{\mu}} 
\begin{cases}
u+\tr^{\Sigma}u-\tr^{\Sigma}u & \text{ on } S_1,\\
\tr^{\Sigma}u+u-\tr^{\Sigma}u &\text{ on } S_2
\end{cases} 
\;
=
\; 
\begin{cases}
u &\text{ on } S_1,\\
u &\text{ on } S_2
\end{cases}
\;
= 
\;
u \in H^1_{\mu}.
$$\\

\noindent \textbf{Step 3: Generalized translation and difference quotient}

\noindent We use the sequence $\Gamma_n$ to introduce a generalized notion of translation, i.e. we show how to construct a shift of the function $u$ by the vector $(h,0,0),$ where $h\in \R.$ To fix the perspective, let us choose $h>0$ to be as in the Step 2. Otherwise, we would translate the function in the opposite direction, which can be clearly done by the same method. Define
\begin{equation}\label{translation}
\Gamma_n^{h,x}(x,y,z):= \Gamma_n(x+h,y,z)=\widetilde{\gamma_n}(x+h,y,z)+\widetilde{\delta_n}(x+h,y,z)-\widetilde{\rho_n}(x+h,y,z).
\end{equation}
Notice that 
$$
\widetilde{\gamma_n}(x+h,y,z)=\widetilde{\alpha_n}(x+h,y,z)-\widetilde{\beta_n}(x+h,y,z)=\alpha_n(x+h,y)-\beta_n(x+h,y),
$$ where 
$$
\alpha_n(x+h,y)\lvert_{S}=
\begin{cases}
\alpha_n(x+h,y) & \text{ on } S_1,\\
\alpha_n(-h,y) &\text{ on } S_2
\end{cases}
\;
=
\;
\begin{cases}
\alpha_n(x+h,y) &\text{ on } S_1,\\
\tr^{\Sigma}[\alpha_n(\cdot_1+h,\cdot_2)] &\text{ on } S_2
\end{cases}
\in H^1_{\mu}.
$$ 
The sequence $\beta_n$ has an analogous representation on the structure $S$. Therefore, the translated sequence $\widetilde{\gamma_n}$ can be expressed as
\begin{align*}
\widetilde{\gamma_n}(x+h,y,z) &=
\begin{cases}
\alpha_n(x+h,y) - \beta_n(x+h,y) &\text{ on } S_1,\\
\alpha_n(-h,y) - \beta_n(-h,y) &\text{ on } S_2 
\end{cases}\\
&=
\begin{cases}
\alpha_n(x+h,y) - \beta_n(x+h,y) &\text{ on } S_1,\\
\tr^{\Sigma}[\alpha_n(\cdot_1+h,\cdot_2)]-\tr^{\Sigma}[\beta_n(\cdot_1+h,\cdot_2)] &\text{ on } S_2.
\end{cases}
\end{align*}

\noindent Passing to the limit in the $H^1_{\mu}$-norm we obtain the following closed-form expression 
\begin{align*}
\widetilde{\gamma_n}(x+h,y,z) \xrightarrow{H^1_{\mu}}
&\begin{cases}
v(x+h,y) - \left(v(x+h,y)-u(x+h,y)\right) &\text{ on } S_1,\\
\tr^{\Sigma}[v(\cdot_1+h,\cdot_2)]-\left(\tr^{\Sigma}[v(\cdot_1+h,\cdot_2)]-\tr^{\Sigma}[u(\cdot_1+h,\cdot_2)]\right) &\text{ on } S_2
\end{cases}\\
=&\begin{cases}
u(x+h,y) &\text{ on } S_1,\\
\tr^{\Sigma}[u(\cdot_1+h,\cdot_2)] &\text{ on } S_2.
\end{cases}
\end{align*}
From now on, to avoid possible confusion, we will use the notation $u_i:=u\lvert_{S_i},$ for $i=1,2$. Thus the last equality can be written as
$$
\widetilde{\gamma_n}(x+h,y,z) \xrightarrow{H^1_{\mu}}
\begin{cases}
u_1(x+h,y) &\text{ on } S_1,\\
\tr^{\Sigma}\left[u(\cdot_1+h,\cdot_2)\right] &\text{ on } S_2.
\end{cases}
$$

Similar calculations allow us to check the convergence of $\widetilde{\gamma_n}$, which was the extension of $u_2 = u\lvert_{S_2}$. Omitting the details, we arrive at 
\begin{align*}
    \widetilde{\delta_n}(x+h,y,z) &= \widetilde{\alpha_n}'(x+h,y,z)-\widetilde{\beta_n}'(x+h,y,z) \\
    &= \alpha_n'(y,z)(x+h)-\beta_n'(y,z)(x+h)\\
    &= \alpha_n'(y,z)(x)-\beta_n'(y,z)(x)
\end{align*}
(recall that $\alpha', \beta'$ are constructed analogously to $\alpha, \beta$, but this time around on $S_2$). Moreover, we have
$$
\widetilde{\delta_n}\lvert_{S}= 
\begin{cases}
\tr^{\Sigma}\left(\alpha_n' - \beta_n'\right) &\text{ on } S_1,\\
\alpha_n'(y,z)- \beta_n'(y,z) &\text{ on } S_2
\end{cases}
\;
=
\;
\begin{cases}
\tr^{\Sigma}\alpha_n' &\text{ on } S_1,\\
\delta_n(y,z) &\text{ on } S_2,
\end{cases}
$$
thus after passing with $n \to \infty$ we obtain
$$
\widetilde{\delta_n}(x+h,y,z)\xrightarrow{H^1_{\mu}}
\begin{cases}
\tr^{\Sigma}u_2 &\text{ on } S_1,\\
u_2(y,z) &\text{ on } S_2.
\end{cases}
$$

Lastly, we consider the sequence $\rho_n$ associated with the extension of $\Sigma$--trace. We have:
$$
\widetilde{\rho_n}(x+h,y,z)=\rho_n(y)(x+h,z)=\rho_n(y)(x,z)=\rho_n(y),
$$
and further
$$
\widetilde{\rho_n}(x+h,y,z) \xrightarrow{H^1_{\mu}}
\begin{cases}
\tr^{\Sigma}u_1 &\text{ on } S_1,\\
\tr^{\Sigma}u_1 &\text{ on } S_2.
\end{cases}
$$

For the convenience of the readers, we summarise all the obtained limits:
\begin{align*}
    \widetilde{\gamma_n}(x+h,y,z) &\xrightarrow{H^1_{\mu}}
    \begin{cases}
        u_1(x+h,y) &\text{ on } S_1,\\
        \tr^{\Sigma}[u_1(\cdot_1+h,\cdot_2)] &\text{ on } S_2,
    \end{cases} \\
    \widetilde{\delta_n}(x+h,y,z) &\xrightarrow{H^1_{\mu}}
    \begin{cases}
        \tr^{\Sigma}u &\text{ on } S_1,\\
        u_2(y,z) &\text{ on } S_2,
    \end{cases} \\
    \widetilde{\rho_n}(x+h,y,z) &\xrightarrow{H^1_{\mu}}
    \begin{cases}
        \tr^{\Sigma}u &\text{ on } S_1,\\
        \tr^{\Sigma}u &\text{ on } S_2.
    \end{cases}
\end{align*}

Recalling the definition of $\Gamma^{h,x}_n$ in \eqref{translation}, we arrive at 
\begin{align*}
\Gamma_n^{h,x} \xrightarrow{H^1_{\mu}} 
&\begin{cases}
u_1(x+h,y)+\tr^{\Sigma}u - \tr^{\Sigma}u &\text{ on } S_1,\\
\tr^{\Sigma}[u_1(\cdot_1+h,\cdot_2)]+u_2(y,z) - \tr^{\Sigma}u &\text{ on } S_2
\end{cases}\\
=
&\begin{cases}
u_1(x+h,y) &\text{ on } S_1,\\
\tr^{\Sigma}[u_1(\cdot_1+h,\cdot_2)]+u_2(y,z) - \tr^{\Sigma}u_1 &\text{ on } S_2
\end{cases}
 \in H^1_{\mu}.
\end{align*}

\noindent This limit will be treated as a generalized translation of $u$. From now on, we denote 
$$
u^{h,x}
:=
\begin{cases}
u_1(x+h,y) &\text{ on } S_1,\\
\tr^{\Sigma}[u_1(\cdot_1+h,\cdot_2)]+u_2(y,z) - \tr^{\Sigma}u_1 &\text{ on } S_2.
\end{cases}.
$$
In a completely analogous way, we introduce the z-direction translation: $u^{h,z}.$

With the notion of translation at hand, we define generalized difference quotients as 
$$
D^{h,x}u:=\frac{1}{h}(u^{h,x}-u), \quad D^{h,z}u:=\frac{1}{h}(u^{h,z}-u).
$$ 
We have 
$$
D^{h,x}u,\; D^{h,z}u \in H^1_{\mu}.
$$ 
Note that $u^{h,x}$ and $u^{h,z}$ coincide with classical translations on $S_1$ and $S_2$, respectively.

\textbf{Step 4:} Higher regularity estimates

In this step, we use the notion of generalized difference quotients to establish higher regularity.

Let $\phi:=-\xi^2D^{-h,x}(u)\in H^1_{\mu}$, where $\xi \in C^{\infty}(\R^3),\; \supp (\xi\lvert_{S}) \subset \subset S$ is the properly chosen cut-off function to be specified later. Due to the Remark \ref{H1 jako testowe}, it is a suitable test function in weak formulation \eqref{weak} and thus in an equivalent formulation \eqref{cieplo3}. Proceeding with the latter one, we have
\begin{align}
\begin{split}
\int_{\Omega}\widetilde{B_{\mu}}\nabla_{\mu}u \cdot \nabla_{\mu}(D^{h,x}\phi) d\mu 
&=\int_{S_1}\widetilde{B_{\mu}}\nabla_{S_1}u \cdot \nabla_{S_1} D^{h,x}\phi dS_1 + \int_{S_2}\widetilde{B_{\mu}}\nabla_{S_2}u \cdot \nabla_{S_2} D^{h,x}\phi dS_2\\
&= \int_{S_1} \widetilde{B_{\mu}}\nabla_{S_1}u \cdot\nabla_{S_1} D^h_x\phi dS_1 + \int_{S_2} \widetilde{B_{\mu}}\nabla_{S_2}u \cdot\nabla_{S_2} D^{h,x}\phi dS_2.
\end{split}
\end{align}
Here $D^h_x$ is a classical difference quotient, which follows from the observation that the generalised translation $\phi^{x,h}$ coincides with a classical one (in $x$ direction still) on $S_1$, as mentioned before. Therefore the integral over $S_1$ already is in a proper form, but we need to show that a constant independent of $h$ uniformly bounds the second term on the right-hand side. We have 
\begin{align}
\begin{split}
&\int_{S_2}\widetilde{B_{\mu}} \nabla_{S_2}u \cdot \nabla_{S_2}D^{h,x}\phi dS_2 =
\int_{S_2} \widetilde{B_{\mu}}\nabla_{S_2}u \cdot \nabla_{S_2}\left[\frac{1}{h}(\phi^{h,x}-\phi)\right]dS_2 \\
= &\int_{S_2} \widetilde{B_{\mu}}\nabla_{S_2}u \cdot \nabla_{S_2} \left[ \frac{1}{h}(\tr^{\Sigma}[\phi_1(\cdot_1+h,\cdot_2)]+\phi_2(y,z)-\tr^{\Sigma}[\phi_1(\cdot_1,\cdot_2)]-\phi_2(y,z))\right]dS_2 \\
= &\int_{S_2} \widetilde{B_{\mu}}\nabla_{S_2}u \cdot \nabla_{S_2} \left[\frac{1}{h}(\tr^{\Sigma}[\phi_1(\cdot_1+h,\cdot_2)]-\tr^{\Sigma}[\phi(\cdot_1,\cdot_2)])\right]dS_2.
\end{split}
\end{align}

We estimate this last integral. In the below computations, the constant $C>0$ varies from line to line and is independent of the parameter $h.$

\begin{align}\label{szacowanko}
\begin{split}
&\left|\int_{S_2} \widetilde{B_{\mu}}\partial_yu \partial_y \left[\frac{1}{h}\left(\tr^{\Sigma}[\phi_1(\cdot_1+h,\cdot_2)]-\tr^{\Sigma}[\phi_1(\cdot_1,\cdot_2)]\right) \right]dydz\right|\\
\leqslant &\left|\int_{S_2} \widetilde{B_{\mu}}\partial^2_yu \left[\frac{1}{h}\left(\tr^{\Sigma}[\phi_1(\cdot_1+h,\cdot_2)]-\tr^{\Sigma}[\phi_1(\cdot_1,\cdot_2)]\right)\right]dydz\right|\\
+&\left|\int_{S_2}\partial_y\widetilde{B_{\mu}} \partial_yu \left[\frac{1}{h} \left(\tr^{\Sigma}[\phi_1(\cdot_1+h,\cdot_2)]-\tr^{\Sigma}[\phi_1(\cdot_1,\cdot_2)]\right)\right]dydz\right|\\
\leqslant \;&C\norm{\partial_yu}_{L^2(S_2)} \left\Vert\frac{1}{h}(\tr^{\Sigma}[\phi_1(\cdot_1+h,\cdot_2)]-\tr^{\Sigma}[\phi_1(\cdot_1,\cdot_2)])\right\Vert_{L^2(S_2)}\\ 
+ &C \norm{\partial_y^2 u}_{L^2(S_2)} \left\Vert\frac{1}{h}(\tr^{\Sigma}[\phi_1(\cdot_1+h,\cdot_2)]-\tr^{\Sigma}[\phi_1(\cdot_1,\cdot_2)])\right\Vert_{L^2(S_2)}\\ 
= &C(\norm{\partial_yu}_{L^2(S_2)}+\norm{\partial_y^2 u}_{L^2(S_2)})\left\Vert\frac{1}{h}(\tr^{\Sigma}[\phi_1(\cdot_1+h,\cdot_2)]-\tr^{\Sigma}[\phi_1(\cdot_1,\cdot_2)])\right\Vert_{L^2(S_2)}\\
\leqslant &C(\norm{\partial_yu}_{L^2(S_2)}+\norm{\partial_y^2 u}_{L^2(S_2)})\left\Vert\frac{1}{h}(\tr^{\Sigma}[\phi_1(\cdot_1+h,\cdot_2)]-\tr^{\Sigma}[\phi_1(\cdot_1,\cdot_2)])\right\Vert_{L^2(\Sigma)}\\
\leqslant\; &C \left(\norm{\partial_y^2 u}_{L^2(S_2)}+\norm{\partial_yu}_{L^2(S_2)}\right) \left\Vert\frac{1}{h}([\phi_1(\cdot_1+h,\cdot_2)]-[\phi_1(\cdot_1,\cdot_2)])\right\Vert_{L^2(S_1)}.
\end{split}
\end{align}
 The second estimate follows from the assumption that $B_{\mu}\lvert_{S_i} \in W^{1,\infty}(S_i),\; i=1,2.$ 
In the second to last inequality, we used the fact that traces are independent of the $z$-variable; it is possible to estimate from above by the integral over $\Sigma.$ The last inequality results from the linearity and the boundedness of the trace operator. This estimate is meaningful as it allows us to change integration over $\Sigma$ to integration over $S_1$ (where the $x$-variable is 'contained'). 
Using the Young inequality with $\eps,$ we obtain
\begin{align}
&C \left(\norm{\partial_y^2 u}_{L^2(S_2)}+\norm{\partial_yu}_{L^2(S_2)}\right) \left\Vert\frac{1}{h}([\phi_1(\cdot_1+h,\cdot_2)]-[\phi_1(\cdot_1,\cdot_2)])\right\Vert_{L^2(S_1)}\\
\leqslant \;&C\left(\frac{1}{\eps^2}(\norm{\partial_y^2 u}_{L^2(S_2)}+\norm{\partial_yu}_{L^2(S_2)})+\eps\norm{D^h_x\phi_1}_{L^2(S_1)}\right).
\end{align}
As in the classical proof, we choose $\eps>0$ small enough and move the term $\eps\norm{D^h_x\phi_1}_{L^2(S_1)}$ to the left-hand side of \eqref{cieplo3}, obtaining the expression of the type
$$\widetilde{C} \norm{D^h_x\phi_1}_{L^2(S_1)} \leqslant C + \int_{\Omega}\widetilde{f} D^{h,x}\phi d\mu,$$
with positive constants $\widetilde{C}, C$ independent of the parameter $h.$

We deal with the right-hand side of the equation analogously.
Firstly, we decompose the integral as
$$
\int_{\Omega}\widetilde{f} D^{h,x}\phi\, d\mu = \int_{S_1}\widetilde{f} D^h_x \phi\, dS_1 + \int_{S_2}\widetilde{f} D^{h,x}\phi \,dS_2.
$$
The integral over $S_1$ is in a standard form and can be treated as in the classical proof. To deal with the integral over $S_2$, we use the analogous estimates as in the computations presented in estimation \eqref{szacowanko}. Applying the standard difference quotients method (i.e. also fixing a function $\xi$) we conclude a uniform boundedness of the first integral; this is clearly explained in the mentioned Section 6.3.1 of the book of Evans \cite{Evans}. In this way, we end up with the formula
$$
\norm{D^h_x \phi_1}_{L^2(S_1)} \leqslant C,
$$ where $C>0$ is again independent of $h.$ 
This implies that
$$
\partial^2_xu \in L^2(S_1)
$$ and $\partial^2_xu$ is a second weak derivative of the function $u.$ We can evoke the analogous procedure to provide extra regularity with respect to the $z$-variable. \\

\pagebreak

Concluding all the reasoning, we arrive at
\begin{tw}\label{wazne}
Assume $\mu \in \mathcal{S}$ is as in \eqref{miary} and $\supp \mu$ is as in \eqref{dyski}. Let $u \in H^1_{\mu}$ be a solution to Problem \eqref{cieplo3} with the coefficients matrix $B_{\mu}\lvert_{S_i}\in W^{1,\infty}(S_i),\; i=1,2.$ Then for $i=1,2$ we have $u\in H^2_{\loc}(S_i).$ $\qedhere$
\end{tw}

The same reasoning can also be applied to establish the analogous result if $\dim S_1 = \dim S_2 = 1$, which we state in the following

\begin{tw}
Let 
\begin{equation*}
S_1:=\{(x,0,0)\in \R^3: x\in [-1,1]\}, \quad S_2:=\{(0,0,z)\in \R^3: z\in [-1,1]\},
\end{equation*}
\begin{equation*}
S:=S_1 \cup S_2, \quad \mu:= \mathcal{H}^1\lfloor_{S_1}+\mathcal{H}^1\lfloor_{S_2}.
\end{equation*}
 Let $u \in H^1_{\mu}$ be a solution to \eqref{cieplo3} on $S$ and $B_{\mu}\lvert_{S_i}\in W^{1,\infty}(S_i),\; i=1,2.$ Then for $i=1,2$ the solution satisfies $u\in H^2_{\loc}(S_i).$
\end{tw}

\begin{dow}
Observe that $S \subset \{(x,0,z)\in \R^3: x,z\in\R\}.$ This allows us to work in a two-dimensional subspace and then trivially extend constructions to $\R^3.$ The proof is based on completely analogous constructions and reasoning as in the two-dimensional regularity theorem proved above. To avoid repetition, we skip the presentation of the proof here.$\qedhere$
\end{dow}
	
It remains to address the case of $\dim S_1 \neq \dim S_2$. However, the method presented in the proof of Theorem \ref{wazne} does not work in this setting. Indeed, let us assume that $\dim S_1 < \dim S_2$; as the components $S_1, \, S_2$ of the low-dimensional structure are transversal (condition b) of Definition \ref{class}), this implies $\dim(S_1 \cap S_2) < \dim S_2 - 1,$ and the $S_1\cap S_2$-trace operator of a function defined on $S_2$ is ill-posed, at least in a classical sense. On the other hand, it turns out that a much simpler extension is available to us.

\begin{prop}\label{zero_extension}
Assume that $\dim S_1 \neq \dim S_2.$ Let $w \in H^1(S_i),\; i=1,2.$ Then the functions 
$$
u:= 
\begin{cases}
w &\text{ on }S_1,\\
0 &\text{ on }S_2
\end{cases}
\quad
\text{and}
\quad 
v:= 
\begin{cases}
0 &\text{ on }S_1,\\
w &\text{ on }S_2
\end{cases}
$$
belong to $H^1_{\mu}.$
\end{prop}
\begin{dow}
Without the loss of generality, we can assume that $\dim S_1 > \dim S_2$. Let us observe that $\capac_2(S_1\cap S_2, S_1)=0$ (the $\capac$ stands for the Sobolev capacity -- for its definition see \cite[Section 4.7]{gar}). This ensures the existence of the sequence $\phi_n \in C^{\infty}(\R^2)$ witnessing this fact.  For the function $w\in H^1(S_1)$ we find the sequence $\alpha_n \in C^{\infty}(\R^2)$ approximating in the $H^1$-norm, and for each term we introduce the $\R^3$-extension $\widetilde{\alpha_n}\in C^{\infty}(\R^3)$ by the formula $\widetilde{\alpha_n}(x,y,z):= \alpha_n(x,y).$ Similarly, we extend the $\phi_n$ sequence to the sequence $\widetilde{\phi_n} \in C^{\infty}(\R^3),\; \widetilde{\phi_n}(x,y,z):=\phi_n(x,y).$ Now we define another sequence $\Psi_n\in C^{\infty}(\R^3),\; \Psi_n:= 1-\widetilde{\phi_n}\widetilde{\alpha_n}.$ It can be easily verified that $\Psi_n \xrightarrow{H^1_{\mu}} u,$ what gives $u \in H^1_{\mu}.$   
This provides that $u \in H^1_{\mu}.$ The fact that $v \in H^1_{\mu}$ is established analogously.  

$\qedhere$
\end{dow}

With the above proposition at hand, we can give a short proof of the higher regularity of weak solutions in the instance of components of different dimensions.

\begin{tw}\label{sobreg}
Assume that $\dim S_1 \neq \dim S_2.$ If $u \in H^1_{\mu}$ is a weak solution of Problem \ref{cieplo3} with $B_{\mu}\lvert_{S_i}\in W^{1,\infty}(S_i),\; i=1,2,$ then $u \in H^2_{\loc}(S_i),\; i=1,2.$
\end{tw}
\begin{dow}
For $i=1,2$ denote $b^i:=-D^{-h}(\xi^2 D^h u^i),$ with a smooth function $\xi$ chosen as in the classical proof. By Proposition \ref{zero_extension}, we know that the functions 
$$
a^1:=
\begin{cases}
b^1, &S_1\\
0, & S_2
\end{cases} 
\quad
\text{and}
\quad
a^2:=
\begin{cases}
0, &S_1\\
b^2, &S_2
\end{cases}
$$ are in $H^1_\mu$ and therefore we can use them as test functions. This implies that the considered low-dimensional issue reduces to the classical one, and the classical proof yields the desired regularity. $\qedhere$
\end{dow}

As a last effort in this section, we address two assumptions previously made to simplify the proof of Theorem \ref{wazne}. Firstly, we assumed that without the loss of generality, we can restrict our attention to structures composed of two sub-manifolds. Secondly, we considered only "straightened-out" domains. In what follows, we show that our construction easily carries on to the general setting.

Let us begin with discussing how to pass from the case of two intersecting manifolds, both being subsets of either $\R$ or $\R^2$, to a general case of a low-dimensional structure consisting of multiple components (which we assume to be "straightened-out" still). Let 
$$
S=\bigcup_{i=1}^m S_i \quad \text{and} \quad \mu:= \sum_{i=1}^m\mathcal{H}^{\dim S_i}\mres_{S_i}.
$$
Keeping in line with previous notation, we will denote for $1 \leqslant i\neq j \leqslant m$ 
$$
\Sigma_{ij}:= S_i\cap S_j.
$$ 
Recall that the intersections $\Sigma_{ij}$ are mutually isolated, i.e. for each two $\Sigma_{ij} \neq \Sigma_{i'j'}$ there exist fixed open (in $\R^3$) sets $O_{ij}$ and $O_{i'j'}$ satisfying 
\begin{equation}\label{pokrywka}
\begin{aligned}
\Sigma_{ij} \subset O_{ij}, \quad \Sigma_{i'j'} \subset O_{i'j'} \quad \text{and} \quad O_{ij} \cap O_{i'j'} = \emptyset.
\end{aligned}
\end{equation} 
Let us put 
$$
C^{\infty}_{ij}:= \{\phi \in C^{\infty}(\R^3): \supp\phi\subset O_{ij}\} \quad \text{and} \quad \mu_{ij}:= \mu\mres_{O_{ij}} \text{ for } 1\leqslant i,j \leqslant m.
$$
Then if $u \in H^1_{\mu}$ solves \eqref{cieplo3} it also satisfies 
\begin{equation}
\int_{O_{ij}}\widetilde{B_{\mu}}\nabla_{\mu_{ij}}u \cdot \nabla_{\mu_{ij}} \phi d\mu_{ij} = \int_{O_{ij}}\widetilde{f} \phi d\mu_{ij}
\end{equation}
for all $\phi \in C^{\infty}_{ij}.$ Now we can conclude that if $S_{ij}:= O_{ij} \cap S$ and we know that $u \in H^2_{\loc}(S_{ij}),$ then $u \in H^2_{\loc}(S_i) $ for all $1\leqslant i \leqslant m$. 

Secondly, we address the "straightened-out" assumption. It was justified since after decomposing the general structure $S$ into substructures $S_{ij},\; 1\leqslant i,j \leqslant m$ it is possible to use the proper composition of diffeomorphisms to obtain the demanded "flat" structures. As we work with closed manifolds, the diffeomorphic changes of manifolds produce some smooth densities that are bounded and isolated from zero. This means that such changes do not impact the convergence in the norms of the spaces $H^1_{\mu},$ $H^2(S_i)$ or $D(A_{\mu})$ (the latter space was introduced in paper \cite{Bouchitte}; for convenience of the readers, we recall these notions in the Appendix). For the construction and detailed discussion of the related diffeomorphisms, see \cite[Theorem 2, paragraph 5]{chom}.

Collecting all the elements of our reasoning, we arrive at the
\begin{tw}\label{globreg}
Let 
$$
S = \bigcup_{i=1}^m S_i \quad \text{and} \quad \mu = \sum_{i=1}^m\mathcal{H}^{\dim S_i}\mres_{S_i}.
$$
If $u \in H^1_{\mu}$ is a weak solution to  \eqref{cieplo3} on $S$ with $B_{\mu}\lvert_{S_i}\in W^{1,\infty}(S_i),\; 1 \leqslant i \leqslant m.$ Then ${u\in H^2_{\loc}(S_i)}$ for $1 \leqslant i \leqslant m.$ $\qedhere$
\end{tw}

\section{Global continuity of weak solutions}

In this section, we show how we can apply the main regularity result to obtain continuity of solutions in the case of structures with a constant dimension of components. In other words, we prove that if a measure $\mu$ belongs to the class $\mathcal{S}$ and moreover
\begin{equation}\label{const_dim}
\supp\mu = S=\bigcup_{i=1}^mS_i, \quad \dim S_1 = \ldots = \dim S_m = k = 1, 2,
\end{equation} then weak solutions to elliptic Problem \ref{cieplo3} are continuous.
 
\begin{tw}(Continuity of solutions)\label{ciag}\\
Let $S=S_1\cup S_2,\; \dim S_1 = \dim S_2.$ Let $u \in H^1_{\mu}$ be a weak solution of Problem \eqref{weak} and ${B_{\mu}\lvert_{S_i}\in W^{1,\infty}(S_i),}$ for  $i = 1, 2.$ Then $u \in C(S).$
\end{tw}
\begin{dow}
By the results of Section 3, we know that $u \in H^2_{\loc}(S_i),\; i=1,2.$ By Proposition \ref{slady}
$$
\tr^{\Sigma}u_1=\tr^{\Sigma}u_2 \quad \text{a.e. on } \Sigma.
$$  As $u\in H^2_{\loc}(S_i)$ implies $u \in C(S_i),$ we have $\tr^{\Sigma}u_i=u_i\lvert_{\Sigma}.$ We conclude that $u_1\lvert_{\Sigma}=u_2\lvert_{\Sigma}$ and thus $u \in C(S).$  $\qedhere$    
\end{dow}

\begin{uw}
We make the following observations regarding the continuity of $u$:
\begin{enumerate}[label=(\alph*)]
\item[a)] In the case $d=1,$ we do not use the fact that $u$ is a weak solution to \eqref{weak}. This shows that continuity is a general property of the space $H^1_{\mu}$ considered on structures built with one-dimensional components.
\item[b)] In the case $d=2$ there exist examples of functions which belong to $H^1_{\mu}$ and are discontinuous. In light of Theorem \ref{ciag}, any possible discontinuities might come only from the lack of continuity of $u|_{S_i} \in H^1(S_i)$ for some component $S_i.$
\item[c)] For a general $\mu \in \mathcal{S}$, not necessarily satisfying $\dim S_1 = \ldots = \dim S_m$, we cannot ensure that $u\in H^1_\mu$ solving \eqref{weak} is continuous. Intuitively, this is due to a fact that intersection of some components will be a point, and a single point has zero $W^{1,2}$-capacity in a plane.
\end{enumerate}
\end{uw}

\noindent Note that Theorem \ref{ciag} addresses the case of structures consisting of two components only. As it turns out, we can easily extend this result.

\begin{tw}\label{ciaglosc_ogolne}
Assume that $S$ and $\mu \in \mathcal{S}$ satisfy \eqref{const_dim} and let $u\in H^1_{\mu}$ satisfy \eqref{weak} with ${B_{\mu}\lvert_{S_i}\in W^{1,\infty}(S_i),}$ $1 \leqslant i \leqslant m.$ Then $u \in C(S)$. 
\begin{dow}
By Theorem \ref{ciag} we have that $u \in C(O_{ij})$ for $1 \leqslant i,j \leqslant m,$ where the sets $O_{ij}$ are elements of the covering as in \eqref{pokrywka}. Regularity Theorem \ref{sobreg} implies that $u \in H^2_{\loc}(S_i),\; 1 \leqslant i \leqslant m,$ thus $u \in C(S_i).$ This two results immediately yield $u \in C(S).$ $\qedhere$
\end{dow}
\end{tw}

\section{Membership in the domain of the second $\mu$-derivative operator}

Having established the higher Sobolev-type regularity, it is only natural to consider whether it is possible to obtain a strong form of the original equation. To this end, it is necessary to introduce second-order differential operators defined on structures in $\mathcal{S}$. It turns out to be quite an involved and delicate matter; for that reason, we skip the presentation here and refer the reader to the Appendix for a brief exposition of the definitions or to the original paper \cite{Bouchitte} of Bouchitt{\'e} and Fragal{\`a} for a detailed discussion. We keep the notation established in  \cite{Bouchitte} and focus on showing that thanks to the regularity result of Section 3, it is possible to show that a weak solution $u$ of the elliptic problem belongs to the domain of the low-dimensional second-order differential operator. 



Firstly, we need the basic notions of the low-dimensional second-order framework.
The definition of the second-order operator $A_{\mu},$ its domain $D(A_{\mu}),$ the operator $\nabla^2_{\mu},$ the Cosserat vector field $b$ and short discussion of their basic properties can be found in the Appendix.


Now, we introduce the operator $L_{\mu}$, which may be interpreted as the right-hand side of the second-order equation.

This section assumes that $\mu \in \mathcal{S}$ satisfies condition \eqref{const_dim}.

\begin{defi}(Operator of second order equation, \cite[Definition 1]{chom})\label{drugie}$ $\\
Let $D(L_{\mu}) = D(\nabla^2_{\mu})$ and assume that a matrix 
$$
B=(b_{ij})_{i,j = 1}^3, \quad b_{ij} \in C^{\infty}(\mathbb{R}^3)
$$ satisfies Definition \ref{relax}. The tangent divergence operator is defined as
$$
{\diverg_{\mu}:(H^1_{\mu})^3 \to L^2_{\mu}}, \quad \diverg_{\mu} ((v_1,v_2,v_3)) := \tr\nabla_{\mu}(v_1,v_2,v_3),
$$ where $\tr$ stands for the trace of a matrix.
Furthermore, we introduce the operator 
$$
L_{\mu} : D(L_{\mu}) \to L^2_{\mu}, \quad L_{\mu}u = \sum_{i,j=1}^3 b_{ij}(\nabla^2_{\mu}u)_{ij} + \sum_{i=1}^3 (\nabla_{\mu}u)_i \diverg_{\mu}(b_{i1},b_{i2},b_{i3}),
$$ which should be understood as a realization of the operator $\diverg(B \nabla_{\mu} u)$ in the setting of low-dimensional structures. Moreover, in the special case of $B=\text{Id},$ where $\text{Id}\in \R^{3\times3}$ is the identity matrix, we denote the operator $L_{\mu}$ as $\Delta_{\mu}$ and its domain as $D(\Delta_{\mu}).$

\qed 
\end{defi}

\noindent Before proceeding, we recall an important property of $L_\mu$.
\begin{tw}\label{domk}(Closedness of $L_{\mu}$, \cite[Theorem 2]{chom})$ $\\
Assume that $\mu \in \mathcal{S}$ satisies condition \eqref{const_dim}. Then the operator $L_{\mu}:D(L_{\mu}) \to L^2_{\mu}$ is a closed operator in the $L^2_{\mu}$-sense.
\end{tw}


The regularity result of Section 3 is local -- we have shown $H^2$-regularity on subsets that are compactly embedded in component manifolds of a low-dimensional structure. To avoid extensive technicalities, we do not aim for up-to-the-boundary regularity and choose to keep working in a local setting. This calls for a further modification. For $u \in H^1_\mu$ being a solution to \eqref{cieplo3}, let $ \intel\overline{S} \subset\subset \intel S$ be chosen in such a way that $\overline{S}$ is still a low-dimensional structure and consider $\overline{\mu}:= \mu\mres_{\overline{S}}$. Abusing notation a bit, we will denote $\overline{\mu}$ as $\mu$. This allows us to write $u \in H^2(S_i),$ where $S_i$ is a component manifold of $\supp \overline{\mu}.$ Let us clearly state that with this procedure, we do not bother with behaviour on the boundary.

As an intermediate step, we need to establish an auxiliary result showing that the smoothness of the force term is propagated to the solution.

\begin{prop}\label{gladkiereg}
Let $u \in H^1_{\mu}$ be a solution to \eqref{weak} with $B_{\mu} = \text{Id}$ and $f\in C^\infty(\R^3)$. Then $u \in C^{\infty}(\R^3).$
\end{prop}
\begin{dow}
Take the covering of the structure $S$ with $\R^3$-open sets $O_{ij}, \; i,j\in \{1,...,m\}$ such that $\Sigma_{ij}=S_i\cap S_j \subset O_{ij}$ and $O_{ij}$ is isolated from $\Sigma_{i'j'}$ if $i'\neq i$ or $j' \neq j.$ 
Let $\phi \in C^{\infty}(\R^3),$  ${\phi_i := \phi\lvert_{S_i} \in C^{\infty}_c(S_i),\; i=1,...,m,}$ and - for simplicity - let us assume that $\supp \phi \cap S \subset O_{12}$. The same argument works for an arbitrary covering element $O_{ij}, \, i\neq j$. As the support of $\phi$ touches only $S_1$ and $S_2,$ we decompose the left-hand side of \eqref{cieplo3} as
$$
\int_{\Omega}\nabla_{\mu}u \cdot \nabla_{\mu} \phi d\mu = \int_{S_1} \nabla_{S_1}u_1 \cdot \nabla_{S_1} \phi_1dS_1 + \int_{S_2} \nabla_{S_2} u_2 \cdot \nabla_{S_2}\phi_2dS_2.
$$ 
After integrating by parts we obtain, for $i=1,2$,
$$
\int_{S_i}\nabla_{S_i} u_i \cdot \nabla_{S_i}\phi_i dS_i = - \int_{S_i}(\Delta_{S_i}u_i)\phi_i dS_i.
$$
Plugging this into equation \eqref{cieplo3} and moving the $S_2$-integral to the right-hand side, we arrive at
\begin{equation}\label{hauhau}
-\int_{S_1}(\Delta_{S_1}u_1)\phi_1dS_1 = \int_{S_2}(\Delta_{S_2}u_2)\phi_2 dS_2 + \int_{S_1}f_1 \phi_1 dS_1 + \int_{S_2}f_2 \phi_2 dS_2.
\end{equation}

\noindent Now, let $\phi_2^\epsilon \in C^\infty(\R^3)$ satisfy
$$
\supp\phi_2^{\eps} \subset \{(x,y,z)\in \R^3: |z| < \eps\} 
\quad 
\text{and} 
\quad 
\phi_2^{\eps}=\phi_2 \text{ on } \left\{(x,y,z)\in \R^3: |z| < \frac{\eps}{2}\right\}
$$
and consider 
$$
\phi^\epsilon := \begin{cases}
    \phi_1 &\text{ on } S_1,\\
    \phi_2^\epsilon &\text{ on } S_2.
\end{cases}
$$
By the construction, we have $\phi^{\eps} \in H^1_{\mu}$, thus it is an accessible test function. Since in equation \eqref{hauhau} no derivative acts on either $\phi_1$ or $\phi_2$, we can easily pass with $\eps$ to zero, obtaining 
$$
\int_{S_2} (\Delta_{S_2}u_2)\phi_2^{\eps}dS_2 \xrightarrow{\eps \to 0}0 
\text{ and } 
\int_{S_2}f_2 \phi_2 dS_2 \xrightarrow{\eps \to 0}0.
$$
As a consequence, this gives 
$$
-\int_{S_1}(\Delta_{S_1}u_1)\phi_1dS_1 = \int_{S_1}f_1\phi_1dS_1.
$$
Since $\phi_1 \in C^{\infty}_c(O_{12} \cap S_1)$ can be chosen arbitrary, we conclude 
$$
-\Delta_{S_1}u_1 = f_1 \text{ a.e. in } O_{12} \cap S_1.
$$ 
This implies that $u_1 \in C^{\infty}(O_{12}\cap S_1).$ We proceed analogously on the second component manifold $S_2$ and obtain $u_2 \in C^{\infty}(O_{12}\cap S_2).$ Substituting $O_{12}$ with other elements of the covering, we conclude $u_i \in C^{\infty}(S_i)$ for all $1 \leqslant i \leqslant m.$ Applying the Whitney Extension Theorem to $u$, we obtain $u \in C^{\infty}(\R^3)$.    

$\qedhere$
\end{dow}

With the above result at hand, we show that $u$ belongs to the domain of $A_\mu$. \\

For clarity in the presentation of the next result, let us abandon the identification of the structure $S = \supp \mu$ with its compactly embedded substructure $\overline{S} = \supp \overline{\mu}.$ 

\begin{tw}\label{nalezenie}
Let $u \in H^1_{\mu}$ satisfy equation \eqref{weak} with $B_{\mu} = \text{Id}.$ Assume that $\overline\mu,\; \supp{\overline\mu}=\overline{S}$ is a low-dimensional measure such that $\intel \overline{S} \subset \subset \intel S$ (in the inherited topology on $S$). Then there exists $b \in L^2_{\mu}(\R^3; T_{\mu}^{\perp})$ such that $(u,b) \in D(A_{\overline\mu}).$
\end{tw}
\begin{dow}
Firstly, we deal with the case of a highly regular right-hand side. Let $f \in \mathring{L^2_{\mu}}$ and additionally let us assume that $f \in C^{\infty}(\R^3).$ By Proposition \ref{gladkiereg} we know that $u\in \mathring{H^1_{\mu}}$ being a solution of 
$$
\int_{\Omega}\nabla_{\mu}u \cdot \nabla_{\mu}\phi d\mu = \int_{\Omega} f\phi d\mu \quad \forall \phi \in C^\infty(\R^3)
$$ has the property that $u\lvert_{\overline{S}}$ can be extended to $u\lvert_{\overline{S}} \;\in C^{\infty}(\R^3),$ where $\overline{S}$ is as before a compactly embedded substructure.   

Now, for a given $h \in \mathring{L^2_{\mu}}$, let $u^h\in \mathring{H^1_{\mu}}$ be a solution to \eqref{weak} with $h$ as the right hand side. Consider a sequence $g_n$ such that
$$
g_n \in C^\infty(\R^3) \cap \mathring{L^2_{\mu}}, \quad g_n \xrightarrow{L^2_\mu} h
$$
and denote solution corresponding to the right-hand side $g_n$ as $v_n \in \mathring{H^1_\mu}.$ In other words, we have 
\begin{align*}
    \int_{\Omega} \nabla_\mu v_n \cdot \nabla_\mu \phi \, d \mu &= \int_{\Omega} g_n \phi \, d\mu, \\
    \int_{\Omega} \nabla_\mu u^h \cdot \nabla_\mu \phi \, d \mu &= \int_{\Omega} h   \phi \, d\mu.
\end{align*}
Subtracting both sides, we obtain 
$$
\int_{\Omega}\nabla_{\mu}(v_n-u^h) \cdot \nabla_{\mu}\phi \,d\mu = \int_{\Omega}(g_n-h)\phi \,d\mu
$$
which yields the estimate
\begin{equation*}\label{nierownosc}
\left|\int_{\Omega}\nabla_{\mu}(v_n-u^h)\cdot \nabla_{\mu}\phi \, d\mu \right| 
\leqslant \int_{\Omega}|g_n-h||\phi| \, d\mu
\end{equation*}
Choosing $\phi:= v_n -u^h \in H^1_\mu$ as a test function, we get
\begin{equation*}\label{nierownosc2}
\int_{\Omega}|\nabla_{\mu}(v_n-u^h)|^2 d\mu\leqslant \int_{\Omega}|g_n-h||v_n-u^h|d\mu
\end{equation*}
and after applying the Young inequality with the $\eps$ we conclude
\begin{equation*}\label{nierownosc3}
\begin{aligned}
C\int_{\Omega}|\nabla_{\mu}(v_n-u^h)|^2 d\mu\leqslant \norm{g_n-h}^2_{L^2_{\mu}}, 
\end{aligned}
\end{equation*}
for some positive constant $C.$ 
Passing with $n$ to infinity, we see that 
$$
\nabla_{\mu} v_n \xrightarrow{L^2_{\mu}}\nabla_{\mu}u^h.
$$
Notice that $\int_{\Omega}v_nd\mu=0$ and $\int_{\Omega}u^hd\mu=0.$ Now, the weak Poincar{\'e} inequality \eqref{weakpoincare} implies that 
$$
v_n \xrightarrow{H^1_\mu} u^h.
$$ 
Now, we restrict further considerations on a compactly embedded $\overline{S} \subset S$ with the corresponding measure $\overline{\mu}.$ 
Theorem \ref{domk} provides that the operator $\Delta_{\overline\mu}: D(A_{\overline\mu}) \to L^2_{\mu}$ is closed. From a definition of the domain $D(A_{\overline\mu})$ (see the Appendix) it follows that $v_n \in C^{\infty}(\R^3)$ implies $v_n \in D(A_{\overline\mu}).$ Since we already know  
$$
v_n\xrightarrow{L^2_{\overline\mu}}u^h, 
\quad 
\Delta_{\overline\mu}v_n = g_n\xrightarrow{L^2_{\overline\mu}}h,
$$ 
we derive 
$$
u^h \in D(\Delta_{\overline\mu}) \text{\; and \;} \Delta_{\overline\mu}u^h=h.
$$ 
By the definition of the domain $D(\Delta_{\overline\mu})$ (see Definition \ref{drugie}), there exists 
$b_{u^h} \in L^2_{\overline\mu}(\R^3; T_{\overline\mu}^{\perp})$ such that $$(u^h,b_{u^h})\in D(A_{\overline\mu}).$$

$\qedhere$
\end{dow}

We decided to introduce the assumption $B_{\mu} = \text{Id}$ in order to simplify the proofs of Proposition \ref{gladkiereg} and Theorem \ref{nalezenie}. Without any qualitative changes, we can follow the methods of those proofs and generalise the results to an arbitrary matrix $B_{\mu}$ satisfying Definition \ref{relax}. For the sake of readability, we omit the details.\\




\section{Appendix}
For the convenience of the readers, we recall here the second-order framework originally introduced in paper \cite{Bouchitte}. Intuitively, the goal is to try and mimic the procedure in which we defined the $\mu$-related gradient $\nabla_{\mu}$ in Section 2 and define the space $M_{\mu}$ playing the role of the second-order counterpart of the tangent bundle $T_{\mu}.$ Now, this construction is more involved as we cannot simply decouple the second-order derivatives from the normal gradients.\\
 
Let $\nabla^2 \phi$ be a matrix of the classical second-order derivatives of a function $\phi \in C^2(\R^3)$ and denote the space of symmetric $3\times3$ matrices as $\R^{3\times 3}_{\text{sym}}$.   
Introduce further the space $\mathcal{Z}_{\mu}$ as
$$
\mathcal{Z}_{\mu}:= \left\{(u,\nabla^{\perp}u,\nabla^2u) \in 
L^2_{\mu} \times L^2_{\mu}\left(\R^3; T_{\mu}^{\perp}\right) \times L^2_{\mu}\left(\R^3;\R^{3\times 3}_{\text{sym}}\right)
\; : \;
u \in C^{\infty}_c(\R^3)\right\}.
$$	
Let $\overline{\mathcal{Z}_{\mu}}$ denote the $L^2_{\mu} \times L^2_{\mu}\left(\R^3; T_{\mu}^{\perp}\right) \times L^2_{\mu}\left(\R^3;\R^{3\times 3}_{\text{sym}}\right)$-closure of  $\mathcal{Z}_{\mu}$ with respect to the inherited product norm. We define the space 
$$
\mathcal{V}_{\mu}:= \left\{z \in L^2_{\mu}(\R^3;\R^{3\times 3}_{\text{sym}}): (0,0,z)\in \overline{\mathcal{Z}_{\mu}}\right\}.
$$ 
Recalling Proposition 3.3 (ii) stated in \cite{Bouchitte} we know that for $\mu$-almost every $x\in \R^3$ there exists a $\mu$-measurable multifunction $M^{\perp}_{\mu}:\R^3 \to P(\R^{3\times 3}_{\text{sym}})$ ($P(X)$ stands for the power set of $X$) satisfying 
$$
\mathcal{V}_{\mu} = \left\{z \in L^2_{\mu}(\R^3;\R^{3\times 3}_{\text{sym}}): z(x) \in M^{\perp}_{\mu}(x) \text{ for } \mu-\text{a.e.}\; x\in \R^3 \right\}.
$$ 
Now the space $M_{\mu}$ is introduced as 
$$
M_{\mu} := \left\{z \in L^2_{\mu}(\R^3;\R^{3\times 3}_{\text{sym}}): z(x) \perp  M^{\perp}_{\mu}(x) \text{ for } \mu-\text{a.e}.\; x\in \R^3\right\},
$$
where $\perp$ means the orthogonality with respect to the standard Euclidean scalar product in $\R^{3\times 3}$.  
 
We do not aim here to discuss further properties of the introduced notions. The interested reader is encouraged to check Lemma 3.2 and Proposition 3.3 in \cite{Bouchitte}.
 
We define the operator $A_\mu$ that appears in Section 5.
Let us introduce the set $D(A_\mu)$ defined as
$$\label{dommu}
D(A_{\mu}):=\left\{(u,b)\in L^2_{\mu} \times L^2_{\mu}(\R^3;T_{\mu}^{\perp}): \exists z\in L^2_{\mu}(\R^3;\R^{3\times 3}_{\text{sym}}) \text{ such that } (u,b,z)\in \overline{\mathcal{Z}_{\mu}}\right\}.
$$
This set plays the role of the domain of the operator $A_\mu.$ 
The operator $A_{\mu}:D(A_{\mu})\to L^2_{\mu}$ is introduced by the formula
$$
A_{\mu}(u,b) := Q_{\mu}(z),
$$
where for $\mu$-a.e. $x\in \R^3,$ the operator $Q_{\mu}(x):\R^{3\times 3}_{\text{sym}} \to M_{\mu}(x)$ is an orthogonal projection onto the space $M_{\mu}(x).$
A normal vector field $b$ such that $(u,b) \in D(A_{\mu})$ is customarily called the Cosserat vector field of $u \in L^2_{\mu}.$
\bigskip

Now we deal with another operator evoked in Section 5 -- the operator $\nabla^2_{\mu}.$
Firstly, we introduce its domain as the space 
$$
D(\nabla^2_{\mu}):= \{u \in H^1_{\mu}: \exists b : (u,b)\in D(A_{\mu})\}.
$$
To construct $\nabla^2_{\mu}$ we need to establish some auxiliary notions first. For $\mu$-a.e. $x \in \R^3$ let ${P_{\mu}^{\perp}(x):\R^3 \to T_{\mu}^{\perp}(x)}$ be an orthogonal projection onto $T_{\mu}^{\perp}(x)$. Furthermore, we introduce the operator $C: D(C) \to L^2_\mu$ acting on the space
$$
D(C) := \{v \in L^2_{\mu}(\R^3;\R^3): P_{\mu}(x)v(x) \in (H^1_{\mu})^3,\; P_{\mu}^{\perp}(x)v(x) \in (H^1_{\mu})^3 \text{ for } \mu-\text{a.e}. \; x\in \R^3\},
$$ 
and given by the formula
$$
Cv := P_{\mu}^{\perp}\nabla_{\mu}(P_{\mu}v) +  P_{\mu}\nabla_{\mu}(P_{\mu}^{\perp}v).
$$ 
With this operator we associate a tensor field $T_C: \mathbb{R}^3 \to \mathbb{R}^{3\times 3}$ defined as
$$
T_C(x)v(x) := (Cv)(x).
$$ 
Finally, the operator $\nabla^2_{\mu}: D(\nabla^2_{\mu}) \to (L^2_{\mu})^{3\times 3}$ is introduced: 
$$
\nabla^2_{\mu}u := P_{\mu}A_{\mu}(u,b)P_{\mu} - T_Cb.
$$ 
Using the Hahn-Banach Theorem, the operator $C: D(C)\to L^2_{\mu}$ can be extended to a continuous operator acting on the space $L^2_{\mu}.$ This ensures that we can act with $T_C$ on any Cosserat field $b$ related to $u \in D(\nabla^2_{\mu})$.
 
Proposition 3.15 of \cite{Bouchitte} states that $\nabla_{\mu}^2$ is independent of a choice of a Cosserat field $b.$ This observation provides well-posedness of the $\nabla_{\mu}^2$ operator.
\bigskip

As our effort, in order to further clarify the role of the operator $A_\mu$ introduced before, we recall  Proposition 3.10 of \cite{Bouchitte}, which describes the $A_\mu$ operator for measures $\mu \in \mathcal{S}$. Firstly, in this case we have 
\begin{equation}\label{prop1}
D(A_{\mu}) \subset \left\{(u,b)\in H^1_{\mu} \times L^2_{\mu}(\R^3;T_{\mu}^{\perp}): \nabla_{\mu}u+b \in H^1_{\mu}\right\}.
\end{equation}
Moreover, the operator $A_{\mu}$ admits a more friendly form of
\begin{equation}\label{prop2}
A_{\mu}(u,b) = \nabla_{\mu}(\nabla_{\mu}u+b).
\end{equation}

\end{document}